\documentclass[reqno,a4paper,12pt]{amsart} 

\usepackage{amsmath,amscd,amsfonts,amssymb}
\usepackage{mathrsfs,dsfont}

\numberwithin{equation}{section}
\numberwithin{figure}{section}

\newcommand\R{\mathbb{R}}

\newcommand\Z{\mathbb{Z}}

\newcommand\gam{\gamma}

\newcommand\lam{\lambda}
\newcommand\Lam{\Lambda}

\newcommand\sig{\sigma}

\newcommand\eps{\varepsilon}

\renewcommand\le{\leqslant}
\renewcommand\ge{\geqslant}
\renewcommand\leq{\leqslant}
\renewcommand\geq{\geqslant}
\newcommand\sbt{\subset}

\newcommand{\ft}[1]{\widehat{#1}}
\newcommand{\dotprod}[2]{\langle #1 , #2 \rangle}

\newcommand{\spec}{\operatorname{spec}}

\newcommand{\cm}{\complement}

\theoremstyle{plain}
\newtheorem{thm}{Theorem}[section]
\newtheorem{lem}[thm]{Lemma}

\newtheorem*{claim*}{Claim}

\newcommand{\thmref}[1]{Theorem~\ref{#1}}
\newcommand{\secref}[1]{Section~\ref{#1}}

\newcommand{\lemref}[1]{Lemma~\ref{#1}}

\theoremstyle{definition}

\newtheorem*{definition*}{Definition}
\newtheorem*{remarks*}{Remarks}
\newtheorem*{remark*}{Remark}

\newenvironment{enumerate-num}
{\begin{enumerate}
\addtolength{\itemsep}{4pt}
}
{\end{enumerate}}

\newenvironment{enumerate-alph}
{\begin{enumerate}
\addtolength{\itemsep}{4pt}
}
{\end{enumerate}}

\begin{document}

\title{Schauder frames of discrete translates in $L^2(\mathbb{R})$}

\author{Nir Lev}
\address{Department of Mathematics, Bar-Ilan University, Ramat-Gan 5290002, Israel}
\email{levnir@math.biu.ac.il}

\author{Anton Tselishchev}
\address{Department of Mathematics, Bar-Ilan University, Ramat-Gan 5290002, Israel}
\address{St. Petersburg Department of Steklov Mathematical Institute, Fontanka 27, St. Petersburg 191023, Russia}
\email{celis\_anton@pdmi.ras.ru}

\date{December 19, 2025}
\thanks{Research supported by ISF Grant No.\ 1044/21.}

\begin{abstract}
We construct a uniformly discrete sequence $\{\lambda_1 < \lambda_2 < \cdots\} \subset \mathbb{R}$ and functions $g$ and $\{g_n^*\}$ in $L^2(\mathbb{R})$, such that every $f \in L^2(\mathbb{R})$ admits a series expansion
\[
f(x) = \sum_{n=1}^{\infty} \langle f, g_n^* \rangle \, g(x-\lambda_n)
\]
convergent in the $L^2(\mathbb{R})$ norm.
\end{abstract}

\maketitle

% =======================================

\section{Introduction}
\label{secA1}

\subsection{}
A system of vectors $\{x_n\}_{n=1}^{\infty}$ in a Banach space $X$
is called a \emph{Schauder basis} if every  $x \in X$  admits a unique series expansion
$x = \sum_{n=1}^{\infty} c_n x_n$ where $\{c_n\}$ are scalars.
 In this case there exist 
 biorthogonal continuous linear functionals $\{x_n^*\}$
 such that  the coefficients of the series expansion are given by
 $c_n = x_n^*(x)$   (see e.g.\  \cite[Section 1.6]{You01}).
 If the series converges unconditionally
 (i.e.\ if it converges for any rearrangement of its terms)
for every  $x \in X$, then $\{x_n\}$ 
is said to be an \emph{unconditional} Schauder basis.

Given a function $g \in L^2(\R)$, we denote its translates by
\begin{equation}
\label{eqTL}
	(T_\lambda g)(x) = g(x-\lambda), \quad \lam \in \R.
\end{equation}
A long-standing problem which is still open,  asks
whether the space $L^2(\R)$ admits
a Schauder basis formed by translates of a single function.
The problem was first posed in the paper \cite{OZ92}, 
where it was also proved that no system of translates
can be an unconditional Schauder basis  in $L^2(\R)$.

A sequence $\Lam = \{\lam_n\}_{n=1}^{\infty}$ of real numbers
is said to be \emph{uniformly discrete} if 
\begin{equation}
\label{eq:ud}
\inf_{n \neq m}  |\lam_m - \lam_n| > 0.
\end{equation}
It was observed in \cite{OZ92} that 
the condition  \eqref{eq:ud} is necessary for 
 a system  of translates 
$\{T_{\lam_n} g \}_{n=1}^{\infty}$
to form a Schauder basis in $L^2(\R)$.

\subsection{}
An example of a uniformly discrete real sequence
 is the sequence of integers $\Z$.
It is well known that if $\Lam = \Z$ then a system
$\{T_\lam g\}_{\lam \in \Lam}$ 
cannot even be \emph{complete} in the space $L^2(\R)$,
see e.g.\ \cite[Example 11.2]{OU16}.
Moreoever, it was conjectured (see \cite[Conjecture 1]{RS95}) 
that the same should be true for any uniformly discrete
sequence $\Lam \sbt \R$. However, this is not the case;
it was proved by Olevskii
\cite{Ole97} that for any ``small perturbation'' of the positive integers,
that is,
\begin{equation}
\label{eq:pertz}
\lam_n = n + \alpha_n, \quad  
0 \ne \alpha_n \to 0 \quad (n \to +\infty)
\end{equation}
there exists  $g \in L^2(\R)$
 such that  the system
  $\{T_{\lam_n} g\}_{n=1}^{\infty}$ 
is complete in $L^2(\R)$.
This shows that the arithmetic structure of the set $\Lam$
plays an important role in the problem.

\subsection{}
If $X$ is a Banach space with dual space $X^*$, then a system
$\{(x_n, x_n^*)\}_{n=1}^\infty$ in $X \times X^*$ is called a 
\emph{Schauder frame} (or a quasi-basis) if every $x\in X$ has a series expansion
	\begin{equation}
		\label{def_Sch_frame}
	x=\sum_{n=1}^\infty x_n^* (x) x_n.
	\end{equation}
If the series \eqref{def_Sch_frame}
converges unconditionally for every  $x \in X$, then 
$\{(x_n, x_n^*)\}$ is called an \emph{unconditional} Schauder frame.
We note that if $\{x_n\}$ is a Schauder basis
with  biorthogonal coefficient functionals $\{x_n^*\}$,
then  $\{(x_n, x_n^*)\}$ is a Schauder frame. 
To the contrary, for a Schauder frame $\{(x_n, x_n^*)\}$  the series expansion 
\eqref{def_Sch_frame} need not be unique  and the coefficient functionals $\{x_n^*\}$ 
 need not be biorthogonal to $\{x_n\}$.	Hence 
	Schauder frames form a wider class of representation systems than Schauder bases.

We proved recently \cite{LT25} that in the space $L^2(\R)$
there do not exist unconditional Schauder frames
consisting of translates, i.e.\   of the form
$\{(T_{\lam_n} g, g_n^*)\}$  where $g \in L^2(\R)$,
$\{\lam_n\}$ is a sequence of real numbers (discrete or non-discrete), 
and $\{g_n^*\}$ are continuous linear functionals on $L^2(\R)$.

On the other hand, in \cite[Section 4]{FPT21} a 
 construction was given of (not unconditional) Schauder frames
 of translates in  $L^2(\R)$. In fact,  it was proved 
 that whenever a system of translates $\{T_\lam g\}_{\lam \in \Lam}$
is complete in $L^2(\R)$,  then there exists a Schauder frame
$\{(T_{\lam_n} g, g_n^*)\}_{n=1}^{\infty}$ 
 such that $\{ \lam_n \} \sbt \Lam$. 
 However the Schauder frames obtained by this construction 
are \emph{highly redundant}, as the sequence
 $\{ \lam_n \}$ is composed
of countably many blocks of finite size,
such that each block gets repeated a higher and higher
 number of times. The  sequence
 $\{ \lam_n \}$ thus ``runs back and forth'' through the set $\Lam$.

\subsection{}
The following question was posed in \cite[Problem 4.4]{OSSZ11}:
does there exist  a Schauder frame in $L^2(\R)$ formed by 
a \emph{uniformly discrete} sequence of translates?
It is clear from the above that the problem cannot be settled by
 combining the results from \cite{Ole97}  and \cite[Section 4]{FPT21}.
 
The main goal of the present paper is nevertheless  to provide 
the above  problem with an affirmative solution.
We will prove the following result:

\begin{thm}
\label{thmA1}
There exist 
a uniformly discrete sequence $\{\lam_1 < \lam_2 < \cdots\} \sbt \R$
 and functions $g$ and $\{g_n^*\}$ in $L^2(\R)$,
 such that  every $f \in L^2(\R)$ admits a series expansion
\begin{equation}
\label{eqC1}
f(x) = \sum_{n=1}^{\infty}   \dotprod{f}{g_n^*} g(x-\lambda_n)
\end{equation}
convergent in the $L^2(\R)$ norm.
\end{thm}

Here $\dotprod{\cdot}{\cdot}$ denotes the usual scalar product
in the space $L^2(\R)$.

The proof   involves techniques 
from the papers \cite{Ole97}, \cite{KO99}, \cite{KO01}, \cite{Ole02}.
It yields, similar to \cite{Ole97}, a sequence $\{\lam_n\}_{n=1}^{\infty}$  
that can be obtained from the positive integers
by a small perturbation.

In fact the   technique allows   to establish some stronger forms of \thmref{thmA1}.
We present these extensions and give some
additional remarks in  \secref{secR1}.

% =======================================

\section{Preliminaries}
\label{secP1}

In this section we fix notation and terminology, and
state some known facts that will be used in the proof
of the main result.

\subsection{}
\label{sec:P2.1}
By a \emph{weight} we mean a function $w \in L^1(\R)$, $w(t) > 0$ a.e.
For a weight $w$ consider the space $L^2_w(\R)$
consisting of all functions $f$ such that
$\int_{\R} |f(t)|^2 w (t) dt < + \infty$.
It is a separable Hilbert space endowed with the scalar product
$\dotprod{f}{g} = \int_{\R} f(t) \overline{g(t)} w(t) dt$.

The following result was proved in \cite{Ole97}.

\begin{lem}[{\cite{Ole97}}]
\label{lemP2.5}
Let $\sig(n) = n + o(1)$, $\sig(n) \neq n$ $(n = 1,2,\dots)$. 
Then for any weight $w_0$ one can find another
 weight $w$, $w(t) \le w_0(t)$ a.e.,  such that 
the exponential system
$\{e^{2 \pi i \sig(n) t}\}_{n=1}^{\infty}$
is complete in the space $L^2_w(\R)$.
\end{lem}

\subsection{}
By a \emph{trigonometric polynomial} we mean a finite sum of the form 
\begin{equation}
\label{eq:P1.1}
P(t) = \sum_{j} a_j e^{2 \pi i \sig_j t} 
\end{equation}
where $\{\sig_j\}$ are distinct real numbers, and $\{a_j\}$
are complex numbers. By the \emph{spectrum} of $P$ we mean the set
$\spec(P) := \{\sig_j : a_j \neq 0\}$.
If we have $\spec(P) \sbt [0, +\infty)$ then 
we say that $P$ is \emph{analytic}.
We use the notation
 $\|\ft{P}\|_q = (\sum_{j}|a_j|^q)^{1/q}$
and  $\|\ft{P}\|_\infty = \max_{j}|a_j|$.

We also denote
$\|P\|_{\infty} = \sup \{|P(t)| : t \in \R\}$,
and observe that
$\|P\|_{\infty} \le \|\ft{P}\|_1$.

The (symmetric) \emph{partial sum}
$S_{r}(P)$ of a trigonometric polynomial
\eqref{eq:P1.1} is defined by 
\begin{equation}
\label{eq:P1.2}
 S_{r}(P)(t) = \sum_{ |\sig_j| \le r} a_j e^{2 \pi i \sig_j t}.
\end{equation}

We observe that if  a trigonometric polynomial $P$
  has integer spectrum, $\spec(P) \sbt \Z$,
  then $P$ is $1$-periodic, that is, 
  $P(t+1) = P(t)$. In this case, we use
$\{ \ft{P}(n)\}$, $n \in \Z$, to denote the Fourier
coefficients of $P$.

The following lemma is taken from \cite[Section 2.2]{Ole02}.

\begin{lem}
\label{lemP2.4}
Given any $\eps>0$ and $q>2$, one can find
an analytic polynomial with integer spectrum
\begin{equation}
P(t) = \sum_{n>0} \ft{P}(n) e^{2 \pi i n t}
\end{equation}
and a measurable set $F \sbt [0,1]$ 
 satisfying the following conditions:
\begin{enumerate-num}
\item \label{tlo:i} $\| \ft{P} \|_{q} < \eps$;
\item \label{tlo:ii} $m([0,1] \setminus F) <  \eps$;
\item \label{tlo:iii} $|P(t) - 1 | < \eps$ on $F$;
\item \label{tlo:iv} Any partial sum
\begin{equation}
S_l(P)(t) = \sum_{n=1}^{l} \ft{P}(n) e^{2 \pi i n t}
\end{equation}
 can be decomposed as a sum   $A_{l} + B_{l}$ of two analytic polynomials
  $A_{l}$, $B_{l}$  with integer spectrum, such that
\begin{equation}
\|A_{l}\|_{L^\infty(F)} < 2, \quad 
\|B_{l}\|_{L^2([0,1])} < \eps.
\end{equation}
\end{enumerate-num}
\end{lem}

For the proof see \cite{KO01}, Lemma 4.1 and Remark 2 on pp.\ 382--383.

\subsection{}
A system of vectors $\{x_n\}$ in a separable Hilbert space $H$ 
is called a \emph{Riesz basis} if it can be obtained as
the image of an orthonormal basis 
under a bounded and invertible linear map.
If $\{x_n\}$ is a Riesz basis then it is
an unconditional Schauder basis, and
the  biorthogonal system
 $\{x_n^*\}$ is also a Riesz basis 
 (see  \cite[Section 1.8]{You01}).

\begin{lem}[{see   \cite[Section 1.10]{You01}}]
\label{lemP2.8}
Let $\{e_n\}$ be an orthonormal basis in a Hilbert space $H$.
If $\{x_n\} \subset H$ satisfies
$\sum_n \|x_n - e_n\|^2 < 1$, 
then $\{x_n\}$ is a Riesz basis
in $H$.
\end{lem}

% =======================================

\section{Proof of \thmref{thmA1}}
\label{secP2}

\subsection{}
Following \cite{Ole97} we begin with a reformulation of the main result
in terms of exponential systems in weighted $L^2_w(\R)$ spaces.

\begin{thm}
\label{thmP2.1}
For any weight $w_0$ we can find 
a weight $w$, $w(t) \le w_0(t)$ a.e.,
a uniformly discrete sequence $\{\lam_1 < \lam_2 < \cdots\} \sbt \R$,
and a sequence of functions $\{e_j^*\} \sbt L^2_w(\R)$,
 such that  every $f \in L^2_w(\R)$ admits a series expansion
\begin{equation}
\label{eq:thmP2.1ser}
f(t) = \sum_{n=1}^{\infty}   \dotprod{f}{e_j^*} \,  e^{2 \pi i \lam_j t}
\end{equation}
convergent in the $L^2_w(\R)$ norm.
\end{thm}

We show that  \thmref{thmP2.1} implies \thmref{thmA1}.
The mapping $U(f) := (f \sqrt{w})^{\wedge}$
is a unitary operator   $L^2_w(\R) \to L^2(\R)$,
where $\ft{\varphi}(x)=\int_{\R} \varphi(t) \exp(-2 \pi i xt) dt$
 is the Fourier transform of $\varphi$.
We take $g := U(1)$ and observe that the exponential
system $\{ e^{2 \pi i \lam_j t}\}$
is mapped by $U$ onto the system of translates $\{T_{\lam_j} g\}$.
Hence if we take $g_j^* := U(e_j^*)$ then we obtain
  \thmref{thmA1} as a consequence of \thmref{thmP2.1}.

\subsection{}
\label{sec:uweight}
We now turn to prove \thmref{thmP2.1}.
Assume that we are given a weight $w_0$.   We use \lemref{lemP2.5}
to choose and fix another weight $u$ satisfying
\begin{equation}
\label{eq:umin}
u(t) \le \min \{ w_0(t), \,
  \tfrac1{10} (1+t^2)^{-1} \}
\end{equation}
and such that there is an exponential system with real frequencies
\begin{equation}
\label{eq:olevcompl}
\{e^{2 \pi i \sig(n) t}\}_{n=1}^{\infty}, \quad 
\sig(n) = n + o(1), \quad 
0 < |\sig(n) - n| < \tfrac1{10}
\end{equation}
which is complete in the space $L^2_u(\R)$.

For a measurable set $E \sbt \R$ we let $m(E)$
be the Lebesgue measure of $E$, and
\begin{equation}
m_u(E) := \int_E u(t) dt.
\end{equation}
Then $m_u$ is a finite measure on $\R$. 
Due to \eqref{eq:umin}, for any $F \sbt [0,1]$ we have
\begin{equation}
\label{eq:mufz}
m_u(F + \Z) \le m(F).
\end{equation}
Similarly,  for any $1$-periodic function $g(t)$ we have
\begin{equation}
\label{eq:intuper}
\|g\|_{L^2_u(\R)} \leq \|g\|_{L^2([0,1])}.
\end{equation}

\subsection{}
The space $L^2_u(\R)$ is a separable Hilbert space. Let us choose some
 orthonormal basis $\{\varphi_k\}_{k=1}^{\infty}$ for the space.
We construct by induction 
a decreasing sequence $\{\eps_k\}$ of positive
numbers, an increasing sequence $\{M_k\}$ of positive numbers, 
and a sequence  of functions $\{\gam_k(t)\}$  satisfying 
\begin{equation}
1 = \gam_0(t) \geq \gam_1(t) \geq \cdots \geq \gam_k(t) > 0,
\quad t \in \R.
\end{equation}
We begin by setting $\eps_0 := 1$ and $M_0 := 1$.
 
At the $k$'th step of the induction, we 
observe that the weighted exponential system 
\begin{equation}
\label{eq:gkmexp}
\{\gam_{k-1}(t) e^{2 \pi i \sig(n) t}\}_{n=1}^{\infty}
\end{equation}
is also complete in the space $L^2_u(\R)$. 
Indeed, suppose that some $f \in L^2_u(\R)$ is orthogonal
to the system \eqref{eq:gkmexp}. Then 
$f \cdot \gam_{k-1}$ is a function in 
$L^2_u(\R)$ orthogonal to the system
\eqref{eq:olevcompl}, which is complete in
$L^2_u(\R)$, so we must have
$f \cdot \gam_{k-1} = 0$ a.e.
In turn this implies that $f = 0$ a.e., which
proves that  also
the system \eqref{eq:gkmexp} 
 is complete in $L^2_u(\R)$.

Hence, given any $\eta_k >0$ we can find a trigonometric polynomial
\begin{equation}
\label{eq:defqk}
Q_k(t) = \sum_{n=1}^{N_k} d_{n,k} e^{2 \pi i \sig(n) t}
\end{equation}
such that
\begin{equation}
\label{eq:phikapprox}
\| \varphi_k  - \gam_{k-1} \cdot  Q_k \|_{L^2_u(\R)} < \eta_k.
\end{equation}

We can then choose $\eps_k > 0$ such that
\begin{equation}
\label{eq:abscont}
 \int_{S}  ( |\varphi_k(t)|^2  + |Q_k(t)|^2)   u(t) dt  < \eta_k^2
\end{equation}
holds whenever $S \sbt \R$ is  a measurable set  satisfying
$m_u(S) <  \eps_k$. By choosing
$\eps_k$ small enough  we can also assume that 
\begin{equation}
\label{eq:epsqk}
\text{(a) \,}
\eps_k < \eta_k^2; \quad  
\text{(b) \,}
\eps_k   \|\ft{Q}_k\|_1   <  \eta_k; \quad 
\text{(c) \,}
M_{k-1} \eps_k < \tfrac1{2} \eps_{k-1}.
\end{equation}

By \lemref{lemP2.4} we can now
 find an analytic   polynomial with integer spectrum
\begin{equation}
\label{eq:pklem}
P_k(t) = \sum_{n>0} \ft{P}_k(n) e^{2 \pi i n t}
\end{equation}
and a set $F_k \sbt [0,1]$ satisfying
\begin{enumerate-num}
\item \label{lo:i} $\| \ft{P}_k \|_{\infty} < \eps_k$;
\item \label{lo:ii} $m([0,1] \setminus F_k) <  \eps_k$;
\item \label{lo:iii} $|P_k(t) - 1 | < \eps_k$ on $F_k$;
\item \label{lo:iv} Any partial sum $S_l(P_k)$
 can be decomposed as a sum   $A_{k,l} + B_{k,l}$ such that
\begin{equation}
\label{eq:aklbkl}
\|A_{k,l}\|_{L^\infty(F_k)} < 2, \quad 
\|B_{k,l}\|_{L^2([0,1])} < \eps_k,
\end{equation}
where  $A_{k,l}$, $B_{k,l}$  are analytic polynomials with integer spectrum.
\end{enumerate-num}

We choose a large positive integer $\nu_k$ (to be specified later)  and set
\begin{equation}
\label{eq:defek}
 E_k = \nu_k^{-1} ( F_k + \Z), 
 \end{equation}
 then using
 \eqref{eq:mufz}  we can deduce from \ref{lo:ii}    that 
\begin{equation}
\label{eq:muekcsmall}
  m_u(E_k^\cm) < \eps_k,
\end{equation}
  while due to \ref{lo:iii} we have 
\begin{equation}
\label{eq:pkvkest}
|P_k(\nu_k t) - 1 | < \eps_k, \quad t \in E_k.
\end{equation}

  We now set $M_k := \max \{ M_{k-1} ,  \|P_k\|^2_\infty \|Q_k\|^2_\infty\}$,
  and choose $\gam_k$ to satisfy
\begin{equation}
\label{eq:gmeq}
  \gam_k(t) = \gam_{k-1}(t), \quad t \in E_k,
\end{equation}
and
\begin{equation}
\label{eq:gmless}
  0 < \gam_k(t) (1 + \|P_k\|_{\infty} + \max_{l} \|A_{k,l}\|_\infty) \leq \gam_{k-1}(t), \quad t \in E_k^\cm,
\end{equation}
which completes the inductive step of the construction.

  \subsection{}
  Let now
\begin{equation}
\label{eq:dkdef}
D_k := \bigcap_{j \geq k} E_j.
\end{equation}
By \eqref{eq:epsqk}(c), \eqref{eq:muekcsmall} we have
\begin{equation}
\label{eq:dkest}
\sum_{k=1}^{\infty}  m_u(E_k^\cm) \le
\sum_{k=1}^{\infty}  \eps_k < + \infty,
\end{equation}
hence the Borel-Cantelli lemma implies that
$\bigcup_{k=1}^{\infty} D_k$ is a set 
of full measure in $\R$. 

Define
\begin{equation}
\label{eq:gaminf}
\gam(t) := \inf_{k} \gam_k(t), \quad t \in \R.
\end{equation}
It follows from \eqref{eq:gmeq} that
\begin{equation}
\label{eq:gamstepk}
\gam(t) = \gam_{k-1}(t) > 0, \quad t \in D_k.
\end{equation}
As a consequence, we conclude that  $\gam(t) > 0$ a.e.

  \subsection{}
Next we claim that the estimate
\begin{equation}
\label{eq:pertphk}
\| \varphi_k - \gam \cdot  \tilde{P}_k \cdot Q_k \|_{L^2_u(\R)} < 4\eta_k,
\quad
\tilde{P}_k(t) := P_k(\nu_k t),
\end{equation}
holds for all $k$.  Indeed, we have
\begin{equation}
 \| \varphi_k - \gam \cdot   \tilde{P}_k \cdot Q_k \|_{L^2_u(\R)}^2  
 = \int_{\R} | \varphi_k - \gam \cdot  \tilde{P_k} \cdot  Q_k |^2 u 
\le  \int_{D_k} + \int_{E_k^\cm} + \int_{D_{k+1}^\cm}.
\end{equation}
We estimate each of the three integrals separately. 

First, due to \eqref{eq:pkvkest},  \eqref{eq:gamstepk}
 we have
\begin{align}
& | \varphi_k(t) - \gam(t)  \tilde{P}_k(t) Q_k(t) |^2  \\
& \qquad  \leq
\big( | \varphi_k(t) - \gam(t)  Q_k(t) | +
|  \gam(t) Q_k(t)  ( 1 - P_k(\nu_k t) ) | \big)^2 \\
& \qquad  \leq
2| \varphi_k(t) - \gam_{k-1}(t)  Q_k(t) |^2 +
2 |Q_k(t)|^2 \eps_k^2 , \quad t \in D_k,
\end{align}
hence using \eqref{eq:phikapprox}, \eqref{eq:epsqk}(b)
and the fact that $\int_{\R} u \le 1$ we obtain
\begin{equation}
\label{eq:J1}
\int_{D_k}  | \cdots |^2 u
\leq 2 \|  \varphi_k - \gam_{k-1} \cdot Q_k \|_{L^2_u(\R)}^2  
 + 2 \|Q_k\|_{\infty}^2 \eps_k^2   \le 4 \eta_k^2.
\end{equation}

Second, due to \eqref{eq:gmless} we have
\begin{align}
&  | \varphi_k(t) - \gam(t)  P_k(\nu_k t) Q_k(t) |^2  \\
& \qquad \le 2  | \varphi_k(t)|^2 + 2 |\gam(t)  Q_k(t) |^2  \|P_k\|^2_{\infty} \\
& \qquad \le 2  | \varphi_k(t)|^2 + 2  | Q_k(t) |^2, \quad t \in E_k^\cm,
\end{align}
so using  \eqref{eq:abscont}, \eqref{eq:muekcsmall} we obtain
\begin{equation}
\label{eq:J2}
\int_{E_k^\cm}  | \cdots |^2 u    \leq   2  \eta_k^2.
\end{equation}

Lastly, we have
\begin{align}
&  | \varphi_k(t) - \gam(t)  P_k(\nu_k t) Q_k(t) |^2  \\
& \qquad \le 2  | \varphi_k(t)|^2 + 2 |\gam(t)  P_k(\nu_k t) Q_k(t)|^2 \\
& \qquad \le 2  | \varphi_k(t)|^2 + 2 M_k.
\end{align}
 Observe that by    \eqref{eq:epsqk}(c)  and  \eqref{eq:muekcsmall}    we have
\begin{equation}
m_u ( D_{k+1}^\cm ) \le
\sum_{j > k} m_u(E_j^\cm) \le
\sum_{j > k} \eps_j \le \eps_{k} M_k^{-1},
\end{equation}
so using \eqref{eq:abscont},  \eqref{eq:epsqk}(a) we obtain
\begin{equation}
\label{eq:J3}
 \int_{D_{k+1}^\cm}  | \cdots |^2 u   \leq 2 \eta_k^2  + 2 \eps_k  \le  4 \eta_k^2.
  \end{equation}

Combining \eqref{eq:J1}, \eqref{eq:J2} and \eqref{eq:J3} we 
arrive at the desired estimate  \eqref{eq:pertphk}.

\subsection{}
Suppose now that we 
 choose the sequence $\{\eta_k\}$ to satisfy
16 $\sum_{k=1}^{\infty} \eta_k^2 < 1$.
 Then \lemref{lemP2.8} and 
 the estimate \eqref{eq:pertphk} imply that
the system 
\begin{equation}
\label{eq:S7}
\{ \gam \cdot  \tilde{P}_k \cdot Q_k  \}_{k=1}^{\infty} 
  \end{equation}
forms a Riesz basis in the space $L^2_u(\R)$.
Hence every $f \in L^2_u(\R)$ has a series expansion 
\begin{equation}
\label{eq:S1}
 f = \sum_{k=1}^{\infty} \dotprod{f}{\psi_k} \, \gam \cdot  \tilde{P}_k \cdot Q_k
  \end{equation}
  where $\{\psi_k\} \sbt L^2_u(\R)$ is  the system biorthogonal 
  to \eqref{eq:S7}   and the  convergence of the series \eqref{eq:S1} is unconditional.
   Moreover, $\{\psi_k\}$ is also 
  a Riesz basis in $L^2_u(\R)$, hence
\begin{equation}
\label{eq:S2}
 \sum_{k=1}^{\infty} |\dotprod{f}{\psi_k}|^2 \le K \|f\|^2
  \end{equation}
  where $K$ is a constant not depending on $f$.

Here and below, by $\dotprod{\cdot}{\cdot}$ and
 $\| \cdot \|$ with no subindex, we mean the scalar product and the  norm in $L^2_u(\R)$.

\subsection{}
Notice that we have
\begin{equation}
\label{eq:lacex1}
 \tilde{P}_k(t)  Q_k(t) =  P_k(\nu_k t) Q_k(t) = 
\sum_{m>0} \ft{P}_k(m) Q_{k,m}(t),
  \end{equation}
where $Q_{k,m}$ are trigonometric polynomials defined by
\begin{equation}
\label{eq:lacex2}
Q_{k,m}(t) := Q_k(t) e^{2 \pi i m \nu_k t} = 
\sum_{n=1}^{N_k} d_{n,k} e^{2 \pi i (\sig(n) + m \nu_k) t}.
  \end{equation} 
If we choose the sequence $\{\nu_k\}$ growing sufficiently fast, 
we can guarantee that the spectra
of the polynomials $Q_{k,m}$ follow each other, in the sense that
\begin{equation}
\label{eq:lacex3}
\max \spec (Q_{k,m_1}) < \min  \spec (Q_{k,m_2}), \quad m_1 < m_2,
  \end{equation}
and moreover, the spectra
of the polynomials $\tilde{P}_{k} \cdot Q_{k}$
also follow each other, that is,
\begin{equation}
\label{eq:lacex4}
\max \spec (\tilde{P}_{k_1} \cdot Q_{k_1})  < \min  \spec (\tilde{P}_{k_2} \cdot Q_{k_2}),
 \quad k_1 < k_2.
  \end{equation}
(This property is referred to as  ``separation of spectra'' in \cite[Section 2.3]{KO01}.)

We now define 
\begin{equation}
\label{eq:lamdefunion}
\Lam   := \bigcup_{k =1}^{\infty}  \spec ( \tilde{P}_{k} \cdot Q_{k}),
  \end{equation}
then each point $\lam \in \Lam$ has a unique representation as
\begin{equation}
\label{eq:lamptform}
\lam   = \sig(n) + m \nu_k, \quad k \ge 1,  \quad m \in \spec(P_k), \quad 1 \le n \le N_k.
\end{equation}
The set $\Lam$ is  uniformly discrete due to \eqref{eq:olevcompl}.
   
\subsection{}
For $\lam \in \Lam$  given by \eqref{eq:lamptform} 
we define a function $h_\lam^* \in L^2_u(\R)$ by
\begin{equation}
\label{eq:coefffunc}
h_\lam^*(t) :=  \overline{d_{n,k}  \ft{P}_k(m)} \, \psi_k(t).
\end{equation}
We will prove that if the elements of $\Lam$ 
are enumerated
as $\{\lam_j\}_{j=1}^{\infty}$
by increasing order, that is, $\lam_{j+1} > \lam_j$,
then each
 $f \in L^2_u(\R)$ admits a series expansion 
\begin{equation}
\label{eq:schfexp}
 f (t) = \sum_{j=1}^{\infty} \dotprod{f}{h_{\lam_j}^*}
 \gam(t) e^{2 \pi i \lam_j t}
\end{equation}
convergent in the $L^2_u(\R)$ norm.

Indeed, using  \eqref{eq:defqk}, \eqref{eq:pklem},
 \eqref{eq:lacex3}, \eqref{eq:lacex4} and \eqref{eq:coefffunc},
it follows that any partial sum of the series \eqref{eq:schfexp}
can be decomposed as $S' + S'' + S'''$ where
\begin{equation}
\label{eq:partial1}
S'(t) = \sum_{s=1}^{k-1} \dotprod{f}{\psi_s} \gam(t) \tilde{P}_s(t) Q_s(t),
\end{equation}
\begin{equation}
\label{eq:partial2}
S''(t) =   \dotprod{f}{\psi_k} \gam(t) Q_k(t) S_{l}(P_k)(\nu_k t),
\end{equation}
and
\begin{equation}
\label{eq:partial3}
S''' (t) =  \dotprod{f}{\psi_k} \gam(t) \ft{P}_k(l+1) S_r(Q_{k,l+1})(t)
\end{equation}
for some $k$, $l$ and $r$.
  
We have $\|f - S'\| = o(1)$ as $k \to \infty$ due to \eqref{eq:S1}.

To estimate $S''$ we recall the property
 \ref{lo:iv} of the polynomial $P_k$, that is,
 the partial sum $S_{l}(P_k)$
 can be decomposed as $ A_{k,l} + B_{k,l}$ so that
 \eqref{eq:aklbkl} holds. Let
\begin{equation}
\tilde{A}_{k,l}(t) = A_{k,l}(\nu_kt), \quad \tilde{B}_{k,l}(t) = B_{k,l}(\nu_kt),
\end{equation}
then from \eqref{eq:aklbkl}, \eqref{eq:defek} and using   
\eqref{eq:intuper} we obtain
\begin{equation}
\label{eq:aklblktilde}
\|\tilde{A}_{k,l}\|_{L^\infty(E_k)} < 2, \quad 
\|\tilde{B}_{k,l}\|_{L^2_u(\R)} < \eps_k.
\end{equation}
We now have
\begin{equation}
\label{eq:sppsplit}
\|S''\| \le |\dotprod{f}{\psi_k}| \cdot (\| \gam \cdot Q_k \cdot \tilde{A}_{k,l} \|
+ \| \gam \cdot Q_k \cdot \tilde{B}_{k,l} \|).
\end{equation}
To estimate the first term we split the norm as the sum of two integrals
\begin{equation}
\| \gam \cdot Q_k \cdot \tilde{A}_{k,l} \|^2 = \int_{\R}
| \gam(t) Q_k (t) \tilde{A}_{k,l}(t) |^2 u(t) dt = \int_{E_k} + \int_{E_k^\cm}.
\end{equation}
Then using \eqref{eq:phikapprox}, \eqref{eq:aklblktilde} we have
\begin{align}
\int_{E_k} &\leq 4 \int_{E_k} | \gam(t) Q_k (t)  |^2 u(t) dt 
\leq 4 \| \gam \cdot Q_k \|^2
\leq 4 \| \gam_{k-1} \cdot Q_k \|^2 \\
& \leq 4 ( \|\varphi_k \| + \eta_k )^2
=  4 ( 1 + \eta_k )^2  < 16, \label{eq:phietaest}
\end{align}
while  using \eqref{eq:gmless} we get
\begin{align}
\int_{E_k^c} &\leq  \int_{E_k^\cm} | \gam_{k-1} (t) Q_k (t)  |^2 u(t) dt 
\leq   \| \gam_{k-1} \cdot Q_k \|^2   < 4 
\end{align}
(the last inequality is established in the same way as   \eqref{eq:phietaest}).
Hence
\begin{equation}
\label{eq:sppfirst}
\| \gam \cdot Q_k \cdot \tilde{A}_{k,l} \|^2  = \int_{E_k} + \int_{E_k^\cm} < 20.
\end{equation}
To estimate the second term in \eqref{eq:sppsplit}
we use \eqref{eq:epsqk}(b), \eqref{eq:aklblktilde} to obtain
\begin{equation}
\label{eq:sppsecond}
\| \gam \cdot Q_k \cdot \tilde{B}_{k,l} \|^2 
\leq \|Q_k\|^2_\infty \|  \tilde{B}_{k,l}  \|^2 
\leq \|Q_k\|^2_\infty  \eps_k^2  < \eta_k^2.
\end{equation}
We conclude from 
\eqref{eq:sppsplit}, \eqref{eq:sppfirst}, \eqref{eq:sppsecond} 
and \eqref{eq:S2} that
\begin{equation}
\|S''\| \le 10 |\dotprod{f}{\psi_k}|  = o(1), \quad k \to \infty.
\end{equation}

Finally we  estimate $S'''$ using the property
 \ref{lo:i} of the polynomial $P_k$
 and \eqref{eq:epsqk}(b),  as
\begin{align}
|S'''(t)| & \le |\psi_k(f)| \cdot  \| \ft{P}_k \|_{\infty} \cdot  \|S_r(Q_{k,l+1})\|_{\infty} \\
& \le  |\psi_k(f)| \cdot \eps_k \cdot  \|\ft{Q}_k\|_1
\le  |\psi_k(f)| \cdot \eta_k.
\end{align}
Hence, using again \eqref{eq:S2} and
the fact that $\int_{\R} u \le 1$, we obtain
\begin{equation}
\|S'''\|   \le  |\psi_k(f)| \cdot \eta_k =  o(1), \quad k \to \infty
\end{equation}
as well. 

We thus conclude that 
\begin{equation}
\|f - (S' + S'' + S''')\| \le \|f - S'\| + 
\|S''\|  + \|S'''\| = o(1), \quad k \to \infty
\end{equation}
which shows that the partial sums of the series
\eqref{eq:schfexp} indeed converge to $f$ in the norm.

\subsection{}
Now we can finish the proof of \thmref{thmP2.1}.
Define $w(t) := u(t) \gam(t)^2$, then $w$ is 
a weight and $w(t) \le w_0(t)$ a.e.\ due to \eqref{eq:umin}.
The mapping $(Uf)(t) := f(t) / \gam(t)$
is a unitary operator   $L^2_u(\R) \to L^2_w(\R)$,
which maps the weighted exponential system 
$\{ \gam(t) e^{2 \pi i \lam_j t}\}_{j=1}^{\infty}$
onto the unweighted system 
$\{ e^{2 \pi i \lam_j t}\}_{j=1}^{\infty}$.
Hence if we take $e_j^* := U(h_{\lam_j}^*)$ then 
we can conclude that
for every $f \in L^2_w(\R)$, the series expansion
\eqref{eq:thmP2.1ser} is valid with the
convergence in the $L^2_w(\R)$ norm.
This completes the proof of \thmref{thmP2.1},
and consequently  \thmref{thmA1} also follows.

% =======================================

\section{Remarks}
\label{secR1}

In this section we present some extensions 
of \thmref{thmA1} and give  additional remarks.

\subsection{}
As in \cite{Ole97}, the function $g$ in \thmref{thmA1}
may be chosen to be   infinitely smooth. This can be achieved
by choosing the weight $u(t)$ in \eqref{eq:umin} to be fast decreasing.
On the other hand, according to
 \cite[Proposition 2.2]{OSSZ11}, if a system of uniformly discrete translates
of a function $g$ forms a Schauder frame in $L^2(\R)$,
then $\int_{\R} |g| = + \infty$, so that $g$ cannot have fast decay.

\subsection{}
It follows from \eqref{eq:olevcompl}, \eqref{eq:lamptform}
that the  sequence $\{\lam_j\}$ in \thmref{thmA1} 
is obtained by a small perturbation of a certain subsequence $\{n_j\}$
of the positive integers. In fact, we can 
have  $\lam_j  =  n_j + o(1)$. This can be deduced from 
the fact that in \lemref{lemP2.5}, we may assume that any subsystem
$\{e^{2 \pi i \sig(n) t}\}_{n > N}$
is still complete in $L^2_w(\R)$
(see also \cite{Lev25}).

On the other hand, note that the perturbations cannot decay arbitrarily fast.
In particular, one cannot have 
$\lam_j = n_j + O(n_j^{-\alpha})$, $\alpha>0$,
see \cite[Section 4.1]{KO03}.

\subsection{}
As in \cite{KO01}, \cite{NO09}, the
  sequence $\Lam = \{\lam_j\}$ in \thmref{thmA1} 
may be chosen to be not only
uniformly discrete, but even lacunary, 
that is, it can satisfy the condition
\begin{equation}
\label{eqR2}
\lam_{j+1} - \lam_j \to +\infty, \quad  j \to +\infty.
\end{equation}
Moreover, $\Lam$ can have any ``subexponential'' growth,
i.e.\ it can satisfy  $\lam_{j+1} / \lam_j > 1 + \eps_j$
where $\{\eps_j\}$ is any pregiven positive sequence 
tending to zero, no matter how slowly.
The latter condition is sharp, that is,
 it cannot be replaced by the Hadamard lacunarity condition 
 $\lam_{j+1} / \lam_j > c > 1$
 (see \cite[Section 11.4.1]{OU16}).

\subsection{}
We can  add to the statement of \thmref{thmA1} the condition
 that the coefficients of the series expansion 
\eqref{eqC1} satisfy for every $q>2$ the inequality
\begin{equation}
\label{eqC3}
\big(\sum_{n=1}^{\infty} |\dotprod{f}{g_n^*}|^q \big)^{1/q} \leq K_q \|f\|,
\end{equation}
where $K_q$ is a constant which depends on $q$
but does not depend on $f$. This strengthens the result
from \cite{NO09}. The condition 
is  sharp, that is,   \eqref{eqC3} cannot hold for  $q=2$.

Indeed, to get condition \eqref{eqC3} it suffices
to fix a sequence $q_k > 2$, $q_k \to 2$, and
 choose the polynomial $P_k$ in \eqref{eq:pklem}
to satisfy
$\|\ft{P}_k\|_{q_k} \|\ft{Q}_k\|_{q_k} < 2^{-k}$
using \lemref{lemP2.4}\ref{tlo:i}.

To the contrary, if \eqref{eqC3} was true with
$q=2$,  it would imply that the restriction mapping
$f \mapsto f|_{[0,1]}$  is a compact operator
$L^2(\R) \to L^2([0,1])$ 
(see the proof of  \cite[Proposition 5.1]{FOSZ14})
which is certainly not the case.
Hence \eqref{eqC3} cannot hold for  $q=2$.

\subsection{}
By the uniform boundedness principle,
 if $\{(T_{\lam_n} g, g_n^*)\}_{n=1}^{\infty}$
 is a Schauder frame in $L^2(\R)$
then the coefficient functionals $\{g_n^*\}$
must be uniformly bounded. We observe that 
the proof of \thmref{thmA1} above, in fact yields a sequence  
$\{g_n^*\}$ such that $\|g_n^*\| \to 0$ as  $n \to +\infty$.
On the other hand, we note that it is also possible to choose
$\{g_n^*\}$ to be \emph{seminormalized}, that is, to satisfy 
$A \le \|g_n^*\| \le B$
where $A,B$ are positive constants not depending on $n$
(while at the same time satisfying   \eqref{eqC3} 
for every $q>2$).

To this end we choose the weight $u(t)$ in \eqref{eq:umin} 
such that $u(t) \le |\chi(t)|^2$, where $\chi$
is the Fourier transform of a nonzero function supported 
on $[-\frac1{10}, \frac1{10}]$.
As in \cite{NO09} this ensures that the exponentials
$\{ e^{2 \pi i \lam_j t}\}_{j=1}^{\infty}$
form a Bessel sequence in $L^2_w(\R)$, 
hence the translates $\{T_{\lam_n} g\}_{n=1}^{\infty}$
form a Bessel sequence in $L^2(\R)$
(see e.g.\ \cite[Section 4.2]{You01} for the definition
of a Bessel sequence and its basic properties).
It then follows by an application of \cite[Lemma 2.4]{BC20}
 (where we take $\{g'_i\}$ to be any orthonormal system, say)
that the sequence  $\{g_n^*\}$ can be replaced with another 
seminormalized sequence, such that
$\{(T_{\lam_n} g, g_n^*)\}_{n=1}^{\infty}$ remains a Schauder frame.
It can be checked that \eqref{eqC3} remains valid as well, 
possibly with a different constant $K_q$.

\subsection{}
The statement of 
\thmref{thmA1} can be strengthened so that
not only the system 
$\{(T_{\lam_n} g, g_n^*)\}_{n=1}^{\infty}$ 
forms a Schauder frame in $L^2(\R)$,
but also the system
$\{(g_n^*, T_{\lam_n} g)\}_{n=1}^{\infty}$  
(where the translates serve as the coefficient functionals)
is at the same time another Schauder frame in $L^2(\R)$.

This can be achieved by choosing the sequences
$\{\eps_k\}$  and $\{\nu_k\}$ with a sufficiently
fast decay and growth respectively.

We note that the Schauder frame property of the system
$\{(g_n^*, T_{\lam_n} g)\}_{n=1}^{\infty}$  
is not an automatic consequence, since the Schauder frame 
$\{(T_{\lam_n} g, g_n^*)\}_{n=1}^{\infty}$ 
is not unconditional
(compare with  \cite[Corollary 3.3]{FOSZ14} and \cite[Theorem 4.2]{LT25}).

\subsection{}
\thmref{thmA1} is true also in higher dimensions,
that is:
\begin{thm}
\label{thmR2}
There is a uniformly discrete sequence $\{\lam_n\}_{n=1}^{\infty}
\sbt  \R^d$, and there are 
 functions $g$ and $\{g_n^*\}$ in $L^2(\R^d)$,
 such that  every $f \in L^2(\R^d)$ admits a series expansion
\eqref{eqC1} where the convergence is in the $L^2(\R^d)$ norm.
\end{thm}

This can be deduced from 
\thmref{thmA1} using cartesian/tensor products.

\subsection{}
We conclude the paper with some remarks on
$L^p(\R)$ spaces. 

It is known
\cite{FOSZ14} that in $L^p(\R)$, $p>2$,
there exist \emph{unconditional} Schauder frames 
of translates, i.e.\ of the form
$\{(T_{\lam_n} g, g_n^*)\}$ where $g \in L^p(\R)$,
$\{\lam_n\} \sbt \R$  and
$\{g_n^*\}$ are continuous linear functionals on $L^p(\R)$.
Moreover, $\{\lam_n\}$ may be chosen to
be an arbitrary unbounded sequence, and in particular,
it may increase arbitrarily fast.

To the contrary, we proved in \cite{LT25} that 
none of the spaces $L^p(\R)$, $1\le p \le 2$, 
admits an unconditional Schauder frame 
consisting of  translates.

In \cite[Section 4]{FPT21},
Schauder frames (not unconditional) 
formed by non-discrete translates
were constructed in $L^p(\R)$, $1 \le p <  \infty$.

The  problem as to whether Schauder frames formed
by a \emph{uniformly discrete} sequence of translates
may exist in the space $L^p(\R)$, $1 < p < 2$, is addressed
in our subsequent paper \cite{LT24}.
\thmref{thmA1}  solves the problem 
affirmatively for the space $L^2(\R)$.
Note that in $L^1(\R)$, a system of uniformly discrete 
translates cannot even be
complete, see e.g.\ \cite[Theorem 1.7]{OSSZ11}.

It is also not known whether the space
$L^p(\R)$, $1 < p < \infty$, admits
a Schauder basis formed by translates of a single function.
(If so, then the translates must be 
uniformly discrete \cite[Theorem 1]{OZ92}).
It  is known  that there are no
unconditional Schauder bases of translates 
in any  of these spaces,  see \cite{FOSZ14}.

% =======================================

\end{document}